\documentclass{amsart}
\usepackage{amsfonts,amssymb,amsmath,amsthm,mathrsfs}
\usepackage{url}
\usepackage[dvips]{epsfig}
\newtheorem{thrm}{Theorem}[section]
\newtheorem{lem}[thrm]{Lemma}
\newtheorem{prop}[thrm]{Proposition}

\newtheorem{remark}[thrm]{Remark}

\newtheorem{definition}[thrm]{Definition}

\begin{document}

\author[C.~A.~Mantica and L.~G.~Molinari]
{Carlo~Alberto~Mantica and Luca~Guido~Molinari}
\address{C.~A.~Mantica: I.I.S. Lagrange, Via L. Modignani 65, 
20161, Milano, Italy -- L.~G.~Molinari (corresponding author): 
Physics Department,
Universit\'a degli Studi di Milano and I.N.F.N. sez. Milano,
Via Celoria 16, 20133 Milano, Italy.}
\email{carloalberto.mantica@libero.it, luca.molinari@mi.infn.it}
\title[Conformally quasi-recurrent manifolds]
{Conformally quasi-recurrent\\ pseudo-Riemannian manifolds}

\begin{abstract}
Conformally quasi-recurrent (CQR)$_n$ pseudo-Riemannian manifolds are investigated, 
and several new results are obtained. It is shown that the Ricci 
tensor and the gradient of the fundamental vector are Weyl compatible 
tensors (the notion was introduced recently by the authors and applies to
significative space-times), (CQR)$_n$ manifolds with concircular fundamental 
vector are quasi-Einstein. For 4-dimensional (CQR)$_4$ Lorentzian manifolds 
the fundamental vector is null and unique up to a scaling, it is an eigenvector of the Ricci 
tensor, and its integral curves are geodesics. Such 
space-times are Petrov type-N with respect to 
the fundamental null vector.
\end{abstract}
\date{18 apr 2014}
\maketitle
\section{\bf Introduction}
Recurrent manifolds were investigated by many geometers: Walker \cite{Walker}
studied manifolds with recurrent Riemann curvature, Adati and Miyazawa 
\cite{Adati} studied conformally recurrent ones. Mc Lenaghan and Leroy 
\cite{McLLe} and then Mc Lenaghan and Thompson \cite{McLTh} investigated 
space-times with complex recurrent conformal curvature tensor; 
they showed that the manifolds are Petrov type D and N and obtained the metric form  
for the case of real fundamental vector. Other results are found in 
\cite{Kaigorodov}, \cite{Khan}  and \cite{Walker}.
Prvanovic introduced the notion of ''quasi-recurrence'' \cite{Prv89}\cite{Prv90}:
\begin{definition}
A $n$-dimensional pseudo-Riemannian manifold is conformally quasi recurrent, 
(CQR)$_n$, if it is not conformally flat, and there is
a non-zero vector field $A_i$ (the fundamental vector) such that the recurrence condition
holds:
\begin{equation}
\nabla_i C_{jklm} = 2 A_i C_{jklm} + A_j C_{iklm} + A_k C_{jilm} + 
A_l C_{jkim} + A_m C_{jkli}.\label{eq1.1}
\end{equation}
\end{definition}
\noindent
$C_{jklm}$ is the Weyl, or conformal, curvature tensor:
\begin{align}
C_{jklm} = R_{jklm} + 
\frac{1}{n-2} (g_{m[j} R_{k]l} + R_{m[j} g_{k]l} ) 
- \frac{R}{(n-1)(n-2)} g_{m[j} g_{k]l} \label{Weyltensor}
\end{align}
$R_{kl} = -R_{mkl}{}^m$ is the Ricci tensor and $R = R^m{}_m$
is the curvature scalar.\\  
The Bianchi identity for the Weyl tensor and the recurrence condition give two simple properties \cite{Prv89}:
\begin{align}
A_m C_{jkl}{}^m=0,  \label{eq2.1}\\
\nabla_m C_{jkl}{}^m =0.  \label{eq2.2}
\end{align}
The contraction of \eqref{eq1.1} with $C^{jklm}$ gives $\nabla_i C^2 = 4A_i C^2$ i.e. the fundamental vector 
$A_i$ is exact if $C^2=C_{jklm}C^{jklm}$ is non-zero \cite{BuchRot}: 
\begin{align}
A_i = \tfrac{1}{2} \nabla_i \log |C|. 
\end{align}
The recurrence condition \eqref{eq1.1} 
arises naturally in the study of conformal 
transformations of the metric \cite{BuchRot}: $\hat g_{kl} = e^{2\sigma} g_{kl}$. They leave 
the Weyl tensor (3,1) unchanged,
$\hat C_{jkl}{}^m =C_{jkl}{}^m$, 
while the Christoffel symbols transform according to:
$\hat\Gamma^m_{ij} =\Gamma^m_{ij} + \delta^m{}_iA_j +\delta^m{}_j A_i - g_{ij} A^m $,
where $A_i =\nabla_i \sigma$. 
The covariant derivative $\hat\nabla_i$ of $\hat C_{jklm} = e^{2\sigma} C_{jklm}$ gives the identity:
\begin{align}
&\hat\nabla_i \hat C_{jklm} - A^p\left [ g_{ij} \hat C_{pklm} +g_{ik} \hat C_{jplm} +g_{il} \hat C_{jkpm} +g_{im} \hat C_{jklp}\right ]\\
&= e^{2\sigma}  \left [
\nabla_i C_{jklm} - 2 A_i C_{jklm} - A_j C_{iklm} - A_k C_{jilm} - A_l C_{jkim} - A_m C_{jkli}
\right ]. \nonumber
\end{align}
If the left-hand-side is zero, the recurrence condition \eqref{eq1.1} is 
precisely obtained. 
The conformal transformation of the metric $\hat g_{ij} = |C| g_{ij}$ maps the recurrence condition to $\nabla_i \hat C_{jklm}=0$. This fact was exploited in \cite{Prv89} (page 203) and \cite{BuchRot} (example 2) to construct explicit CQR metrics.

Manifolds with the condition \eqref{eq1.1} were named  {\em pseudo conformally symmetric}, (PCS)$_n$, by De and Biswas \cite{De1}. 
De and  Mazumdar \cite{De2} proved the following: 
if a (PCS)$_n$ manifold admits a conformal motion, i.e. there are a vector field 
$\xi_j$  and a scalar field $\sigma $ such that $\nabla_i \xi_j +\nabla_j \xi_i = 2\sigma g_{ij} $ (\cite{Stephanietal} page 564), then it is either conformally flat or $\nabla_j\sigma $ is a null vector. In particular,  a (PCS)$_4$ space-time with proper conformal motion is type N or O.\\

In this paper we present several new results for (CQR)$_n$ manifolds.
In Section 2,
we prove that the Ricci tensor and the tensor $\nabla_i A_j$ of (CQR)$_n$
manifolds are Weyl compatible. The notion of Weyl compatibility was 
recently introduced by us in \cite{EXTDS}, \cite{RCT} and \cite{WCT}
and is a property of several space-times of interest in physics
\cite{Deszcz1}\cite{Deszcz2}. We obtain analogous recurrence 
properties for the tensors $C_{ijab}C^{ijcd}$ and 
$C_{ijka}C^{ijkb}$; in particular
we show that the latter can be rescaled to a Codazzi tensor.\\
In Section 3 we prove that (CQR)$_n$ manifolds with concircular fundamental vector
are quasi-Einstein, and are
Ricci pseudo-symmetric in the sense of Deszcz.\\
In Section 4 we investigate a collateral class of pseudo-Riemannian 
manifolds, $A_i C_{jklm} + A_j C_{kilm} + A_k C_{ijlm} = 0$, and show
that all Pontryagin forms vanish. The result is useful for the
discussion in Section 5, devoted to 4-dimensional (CQR)$_4$ Lorentzian 
manifolds (space-times). We prove that the fundamental vector $A$ is null 
and unique up to a scaling, it is an eigenvector of the Ricci 
tensor, and its integral curves are geodesics. Finally we show that such 
space-times are Petrov type-N with respect to $A$.\\

Hereafter, the manifolds are Hausdorff, connected, of dimension
$n\ge 3$, with a Levi-Civita connection ($\nabla_j g_{kl} = 0$), and
non-degenerate metric of arbitrary signature (pseudo-Riemannian 
manifolds. For Lorentzian manifolds the signature is 2).  


\section{\bf (CQR)$_{\mathbf n}$ manifolds: general properties}
The property $\nabla_m C_{jkl}{}^m = 0$ has an interesting consequence. 
In \cite{SecOrd} and \cite{RCT} we proved a differential identity that extends an identity by Lovelock for the Riemann tensor
to curvature tensors; for the Weyl tensor it is:
\begin{lem}
In a pseudo-Riemannian manifold:
\begin{align}
&\nabla_i \nabla_m C_{jkl}{}^m + \nabla_j \nabla_m C_{kil}{}^m + 
\nabla_k \nabla_m C_{ijl}{}^m \label{eq2.6} \\
&= -\frac{n-3}{n-2} (R_{im} R_{jkl}{}^m + 
R_{jm} R_{kil}{}^m + R_{km} R_{ijl}{}^m). \nonumber
\end{align}
\end{lem}
By eq.\eqref{eq2.2} the left-hand-side is zero, and a result that was first proven 
in \cite{Prv89} is here obtained straightforwardly:
\begin{prop}\label{2.4}
In a (CQR)$_n$ pseudo-Riemannian manifold, the Ricci tensor
is Riemann compatible:
\begin{equation}
R_{im} R_{jkl}{}^m + R_{jm} R_{kil}{}^m + R_{km} R_{ijl}{}^m = 0. \label{eq2.7}
\end{equation}
\end{prop}
\begin{remark}
The property of compatibility was introduced in \cite{EXTDS}
and investigated in detail in \cite{RCT} and \cite{WCT}; in the latter it is 
proven that a Riemann compatible tensor is Weyl compatible. 
Then eq.\eqref{eq2.7} implies: 
\begin{equation}
R_{im} C_{jkl}{}^m + R_{jm} C_{kil}{}^m + R_{km} C_{ijl}{}^m = 0.\label{eq2.8}
\end{equation}
It has been shown recently that several important 4-dimensional 
space-times (as Robertson-Walker, or G\"odel metrics) have the property
that the Ricci tensor is Weyl-compatible \cite{Deszcz1}\cite{Deszcz2}.
This also occurs in certain differential structures for the tensor $Z_{ij}=R_{ij}+\varphi g_{ij}$,
that are investigated in \cite{Mant_Mol}\cite{Mant_Suh1}\cite{Mant_Suh2}\cite{Mant_Suh3}.\\
If the manifold has a matter content, the stress-energy 
tensor $T_{kl}$ is linked to the Ricci tensor by Einstein's equation \cite{DeFCl}\cite{Stephani}: 
$R_{kl} -\frac{1}{2}R \,g_{kl} = {\rm G}\, T_{kl}$.
Therefore, eq.s \eqref{eq2.7} and \eqref{eq2.8} imply that in (CQR)$_n$ manifolds
the stress-energy tensor is Riemann and Weyl compatible.\\

\end{remark}

The following statements are about the fundamental vector $A_i$ and its gradient:
\begin{prop}
In a (CQR)$_n$ manifold, with fundamental vector $A_i$, 
\begin{align}
&(\nabla_i A^m) C_{jklm} + (A^m A_m )C_{jkli} = 0; \label{eq2.3} \\
&(\nabla_i C_{jklm})(\nabla^i C^{jklm} ) = 8(A_i A^i ) (C_{jklm}C^{jklm} );
\label{eq2.12}
\end{align}
the tensor $\nabla_i A_m$ is  Weyl compatible:
\begin{equation}
 (\nabla_i A_m )C_{jkl}{}^m + (\nabla_j A_m )C_{kil}{}^m + 
(\nabla_k A_m ) C_{ijl}{}^m = 0. \label{eq2.11W}
\end{equation}
\begin{proof}
%
The first identity was obtained in \cite{Prv89}
and reobtained in \cite{BuchRot}. 
The second one is proven by squaring \eqref{eq1.1}; several terms vanish because of \eqref{eq2.1}. Finally, write three versions of eq.\eqref{eq2.3} with indices $ijk$ cyclically 
permuted and sum up, a cancellation occurs by the Bianchi identity of
the Weyl tensor. 
\end{proof}
\end{prop}
%
%
We now obtain some results for the tensor
$\Theta_{pqrs} = C_{pqlm}C_{rs}{}^{lm}$,  and its contractions
$\Gamma_{pr} = C_{pqlm} C_r{}^{qlm}$ and $\Gamma^k{}_k=C^2$. Besides the 
general symmetries $\Theta_{pqrs}= -\, \Theta_{qprs}= -\, \Theta_{pqsr}=\Theta_{rspq}$ and $\Gamma_{pq}=\Gamma_{qp}$,
it is: $A^p\Theta_{pqrs}=0$ and $A^p \Gamma_{pr}=0$. The recurrence 
property \eqref{eq1.1} of the Weyl tensor gives: 
\begin{align}
\nabla_i\Theta_{pqrs} = 4A_i\Theta_{pqrs} + A_p \Theta_{iqrs} +A_q \Theta_{pirs} +A_r \Theta_{pqis}
+A_s\Theta_{pqri}.
\end{align}
and the contracted ones: 
\begin{align}
&\nabla_s \Theta_{pqr}{}^s = A_q \Gamma_{pr}-A_p \Gamma_{qr}\label{TGG}\\
&\nabla_i \Gamma_{pr} = 4A_i\Gamma_{pr}+A_p\Gamma_{ir}+A_r\Gamma_{ip}. \label{nablagamma}
\end{align}
\begin{prop}
If $\Gamma^q{}_q\neq 0$ then 
$|\Gamma |^{-3/4}\Gamma_{pq}$ is a Codazzi tensor.
\begin{proof} From \eqref{nablagamma} it follows that
\begin{align}
\nabla_j \Gamma_{kl}-\nabla_k \Gamma_{jl} =3A_j\Gamma_{kl} -3A_k \Gamma_{jl}
\label{coda}
\end{align} 
Since $A_i = \frac{1}{4}\nabla_i \log |\Gamma |$, it is 
$\nabla_i(|\Gamma |^{-3/4}\Gamma_{jq})=\nabla_j (|\Gamma |^{-3/4} \Gamma_{iq})$.
\end{proof}
\end{prop}
\noindent
In \cite{EXTDS} we proved that Codazzi tensors are Riemann compatible. Since the
compatibility relation is algebraic, the normalization factor is irrelevant and it follows 
that $\Gamma_{ij}$ is Riemann and Weyl compatible, and it commutes with the Ricci
tensor \cite{RCT}.


Eqs.\eqref{TGG} and \eqref{coda} give an identity analogous
to the contracted second Bianchi identity of the Riemann tensor:
$$\nabla_m \Theta_{jkl}{}^m =\frac{1}{3}(\nabla_k \Gamma_{jl}-\nabla_j\Gamma_{kl} ).$$
Another covariant derivative $\nabla_i$ and summation on cyclic permutations of the
indices $ijk$ produces the property ($C^2\neq 0$):
\begin{align}
\nabla_i\nabla_m\Theta_{jkl}{}^m +
\nabla_j\nabla_m\Theta_{kil}{}^m +\nabla_k\nabla_m\Theta_{ijl}{}^m =0 \label{Love}
\end{align}
%
because $\Gamma_{ij}$ is Riemann compatible. Contraction with $g^{il}$ gives  
$\nabla_l\nabla_m \Theta_{jk}{}^{lm}=0$, as it always occurs for the Riemann tensor.
{\quad}\\

\section{\bf (CQR)$_\mathbf{n}$ with concircular fundamental vector}
A connection of (CQR)$_n$ manifolds with pseudo-symmetric
manifolds of Deszcz-type is shown, if the fundamental vector is concircular. 
We refer to the following definition (more general ones are possible):
\begin{definition}
A vector field $V$ is concircular if $\nabla_s V_i = V_s V_i + \gamma \, g_{si}$,
where $\gamma $ is a scalar field. \end{definition}

\begin{prop} In a (CQR)$_n$ manifold, if the fundamental vector is concircular, then:
$\gamma =-A_mA^m$ and $\gamma $ is constant.
\begin{proof}
If $A$ is concircular, it is $(\nabla_s A^m)C_{jklm} = A_s A^mC_{jklm}+\gamma C_{jkls}$; 
the properties \eqref{eq2.3} and \eqref{eq2.1} imply 
$A^mA_m=-\gamma $. 
Finally $\nabla_i \gamma =-2A_m\nabla_i
A^m = -2A_m(A_iA^m+\gamma \delta_i{}^m) = 0$, i.e. $\gamma $ is a constant.
Let us also note that $\nabla_sA^s =(n-1)\gamma $.
\end{proof}
\end{prop}

We begin with the case $\gamma \neq 0$.
The derivative of the concircularity condition and the Ricci identity give the integrability condition $A_m R_{jkl}{}^m =\gamma (A_k g_{jl} -A_jg_{kl} )$,  then $A_mR_j{}^m =(n-1)\gamma A_j$. If this relation is used in the expression for the Weyl tensor \eqref{Weyltensor},
bearing in mind that $A^m C_{jklm} = 0$, we obtain: 
\begin{align}
R_{kl}= \frac{A_kA_l}{A_mA^m} \left(  n\gamma - \frac{R}{n-1}\right ) + g_{kl} \left ( \frac{R}{n-1}-\gamma \right).
\end{align}
\begin{prop}
A (CQR)$_n$ with non-null concircular fundamental vector is quasi-Einstein.
\end{prop}
Quasi-Einstein Riemannian manifolds were investigated by Chaki and Maity \cite{ChakiMaity}; 
quasi-Einstein pseudo-Riemannian ones occur
in the study of exact solutions of the Einstein equation, and in the study of 
quasi-umbilical surfaces in semi-Euclidean spaces \cite{DeszczGlo}.\\

The case $\gamma = 0$ is now considered; it always occurs in $n=4$. Multiplication by $A^m$ of
\eqref{Weyltensor} and use of the relation $A_m R_{jkl}{}^m =0$,  that implies $A_mR_j{}^m =0$  
and $A_mC_{jkl}{}^m=0$, give:
\begin{align}
A_k \left[ g_{jl} \frac{R}{n-1} -R_{jl} \right ] = A_j \left[ g_{kl} \frac{R}{n-1} -R_{kl} \right]. \label{A=A}
\end{align}
\begin{thrm}
In a (CQR)$_n$, $n \ge 4$, if the fundamental vector is such that $A_iA^i=0$ and 
$\nabla_i A_j = A_i A_j$, then the scalar curvature is zero, and the Ricci tensor is a Codazzi tensor, with rank
not greater than 1. 
\begin{proof}
Multiplication of \eqref{A=A} by $A^j$ gives $R = 0$. 
This means that $0=\nabla_m C_{jkl}{}^m=\nabla_m R_{jkl}{}^m$. Then 
$\nabla_m R_{ij}=\nabla_i R_{mj}$ (the Ricci tensor is a Codazzi tensor).\\  
As $R=0$ it is: $A_k R_{jl} = A_j R_{kl} $.
Let $ \theta^i A_i =1$ then $R_{kl} =A_l \theta^j R_{jk}$ and, by symmetry, 
$A_l \theta^j R_{jk}=A_k \theta^j R_{jl}$ then $\theta^j R_{kj} = 
A_k (\theta^j\theta^l R_{jl} )$ and $R_{kl} =A_kA_l (\theta^j\theta^m R_{jm} )$ i.e. the
Ricci tensor has rank one (if non-zero).
\end{proof}
\end{thrm}

For a (CQR)$_n$ the commutator of two covariant derivatives of the Weyl tensor is
\begin{align*}
(\nabla_i \nabla_s & -\nabla_s\nabla_i) C_{jklm} \\
=\,  2\,&(\nabla_s A_i - \nabla_i A_s)C_{jklm} + (\nabla_s A_j - A_s A_j)C_{iklm} + (\nabla_s A_k - A_s A_k )C_{jilm}\\
+  &(\nabla_s A_l-A_sA_l)C_{jkim} +(\nabla_s A_m - A_sA_m )C_{jkli} - (\nabla_i A_j -A_iA_j)C_{sklm}\\ 
- &(\nabla_i A_k - A_i A_k )C_{jslm} - (\nabla_i A_l - A_i A_l )C_{jksm} - (\nabla_i A_m - A_i A_m )C_{jkls}.
\end{align*}
If the fundamental vector is concircular the commutator becomes the 
condition for the manifold to be Weyl-pseudo-symmetric according to Deszcz \cite{Deszcz1989}: 
\begin{align}
(\nabla_i\nabla_s -\nabla_s\nabla_i)C_{jklm} 
=& \,\gamma \,[ g_{sj}C_{iklm} - g_{ij}C_{sklm} + g_{sk}C_{jilm} - g_{ik} C_{jslm} \label{PSD}\\
& +g_{ls} C_{jkim} - g_{li} C_{jksm} + g_{ms} C_{jkli} - g_{mi} C_{jkls}  ].\nonumber
\end{align}
In the same paper it is proven that for $n\ge 5$ the Weyl tensor in \eqref{PSD} can be replaced with
the Riemann tensor, i.e. the manifold is also (Riemann) pseudo-symmetric according to Deszcz \cite{Deszcz1987}.  Finally, a contraction gives:
$$(\nabla_i\nabla_s -\nabla_s\nabla_i)R_{kl} =\,\gamma\, \left [ g_{sk}R_{il} - g_{ik} R_{sl} +g_{ls} R_{ki} - g_{li} R_{sk}
\right ],$$
i.e. the manifold is Ricci-pseudo-symmetric according to Deszcz \cite{Deszcz89}. \\
In analogy to this, the evaluation of
\begin{align*}
(\nabla_s\nabla_i-\nabla_i\nabla_s)\Gamma_{pr} =&\,
4\, (\nabla_s A_i-\nabla_i A_s)\Gamma_{pr} + (\nabla_s A_p-A_p A_s)\Gamma_{ri}
+ (\nabla_s A_r-A_s A_r)\Gamma_{ip}\\
& - (\nabla_i A_p-A_i A_p)\Gamma_{rs} - (\nabla_i A_r-A_i A_r)\Gamma_{sp}
\end{align*}
simplifies if the fundamental vector is concircular, and the Gamma tensor gains a Deszcz-type recurrent symmetry
\begin{align*}
(\nabla_s\nabla_i-\nabla_i\nabla_s)\Gamma_{pr} =\,\gamma \, [g_{is}\Gamma_{pr} + g_{sp}\Gamma_{ri}
+ g_{sr}\Gamma_{ip} - g_{ip}\Gamma_{rs} - g_{ir}\Gamma_{sp}]
\end{align*}


\section{\bf A related class of manifolds}
We turn for a while to another class of manifolds that have been investigated 
in the geometric literature (see for example \cite{DefDes}\cite{DesHot}\cite{DesGry}):
$n$-dimensional pseudo-Riemannian manifolds with a vector field $A_i$ such that
\begin{equation}
 A_i C_{jklm} + A_j C_{kilm} + A_k C_{ijlm} = 0. \label{eq2.13} 
\end{equation}
The results that are obtained will be useful in the study of (CQR)$_4$ manifolds.\\

The condition allows an explicit representation of the Weyl tensor, first obtained in 
\cite{DefDes}. We prove it for completeness, and add the new result of
Weyl compatibility: 
\begin{prop}\label{prop2.8}
In a non-conformally flat pseudo-Riemannian manifold, 
if \eqref{eq2.13} holds, then: 1)  $A_i A^i = 0$, 2) there is a symmetric tensor 
$E_{kl}$ such that: 
\begin{align}
C_{jklm}=A_jA_m E_{kl}-A_jA_lE_{mk}-A_kA_m E_{jl}+A_kA_lE_{jm} \label{CAAE}
\end{align}
3) the tensor $E_{kl}$ is Weyl compatible:
\begin{align}
 E^m{}_i C_{jklm} + E^m{}_j C_{kilm} + E^m{}_k C_{ijlm} =0. \label{EWEYLC}
 \end{align}
\begin{proof}
1) By contracting \eqref{eq2.13} with $g^{im}$ one obtains $A^mC_{jklm}=0$; 
contraction of \eqref{eq2.13} with $A^i$ gives $(A^i A_i)C_{jklm} = 0$, which
implies $A^iA_i=0$.  \\
2) Let $B^i$ be a unit (or null) vector such that $B^iA_i=1$. A first 
contraction of  \eqref{eq2.13} 
with $B^i$ gives 
$C_{jklm}=A_j B^i C_{iklm} - A_k B^i C_{ijlm} $ (*); another contraction gives 
$B^m C_{mkj}=A_jE_{kl}-A_k E_{jl}$, with $E_{kl}=E_{lk}=B^i B^m C_{iklm}$.
Inserting this back in equation (*) gives the result.\\
3) The representation \eqref{CAAE} gives $E^m_i C_{jklm} =
A_kA_lE^m{}_iE_{jm} -A_jA_l E^m{}_iE_{mk}$. The sum of three versions of it with indices  $ijk$
cyclically permuted gives \eqref{EWEYLC}.
\end{proof}
\end{prop}
Note the properties 
$E^k{}_k=0$, $A^kE_{kl}=0$, $B^k E_{kl}=0$. If $B^k$ is space-like or time-like,
$E_{kl}$ is the {\em electric part} of the Weyl tensor \cite{Bert}.
\begin{prop}\label{prop2.8bis}
In a $n$-dimensional non-conformally flat pseudo-Riemannian manifold, 
if \eqref{eq2.13} holds, then:\\
1) $C_{lmj}{}^k C_{pqk}{}^j = 0$,\\
2) $C_{lma}{}^b C_{pqb}{}^c C_{rsc}{}^d C_{tud}{}^a=0$.
\begin{proof}
Contraction of \eqref{eq2.13} with $C^{kj}{}_{pq}$ and  
the relation $A_m C_{jkl}{}^m = 0$ give  $A_i C_{jklm} $ $C^{kj}{}_{pq} = 0$, 
from which 1) follows. The next statement and higher order ones are proven analogously.
\end{proof}
\end{prop}
The proposition 
is relevant in the study of Pontryagin forms. 
Given orthonormal tangent vectors $X_r$, Pontryagin forms $p_k$ are 
antisymmetric combinations 
$ p_k (X_1 ... X_{4k}) = \sum_P (-1)^P \omega_k ( X_{i_1} ... X _{i_{4k}})$
of forms $\omega_k$  (see \cite{DerdzRot} \cite{DerdzSh} \cite{RCT} \cite{Nakahara} 
\cite{Postnikov}):
\begin{align*}
&\omega_1 ( X_1\ldots X_4 ) = R_{ija}{}^b R_{klb}{}^a ( X^i_1 \wedge X^j_2)
( X^k_3 \wedge X^l_4), \label{eq2.15}\\
&\omega_2 (X_1\ldots X_8) = R_{ija}{}^b R_{klb}{}^c R_{mnc}{}^d R_{pqd}{}^a  
(X^i_1 \wedge X^j_2)\ldots (X^p_7 \wedge X^q_8), \quad \text{etc.} \nonumber 
\end{align*}
It was shown by Avez \cite{Avez} (see also \cite{DerdzRot}) that forms are
unchanged if Riemann's tensor is replaced by Weyl's tensor, for example:
$$ \omega_1 ( X_1.... X_4 ) = C_{ija}{}^b C_{klb}{}^a ( X^i_1 \wedge X^j_2)
( X^k_3 \wedge X^l_4 ).  $$

\begin{prop} \label{prop2.9}
In a pseudo-Riemannian manifold of type \eqref{eq2.13}, all 
Pontryagin forms vanish.
\begin{proof}
By prop.\ref{prop2.8} it is $C_{lmj}{}^k C_{pqk}{}^j = 0$, then $\omega_1=0$
and $p_1=0$. Moreover, by the same proposition it is $C_{lma}{}^bC_{pqb}{}^c C_{rsc}{}^d
C_{tud}^a=0$, then $\omega_2=0$ and $p_2=0$. A trick analogous to that used in 4) of
prop.\ref{prop2.8} shows that all Pontryagin forms vanish.
\end{proof}
\end{prop}

\begin{remark}
For a compact orientable 4-dimensional pseudo-Riemannian manifold 
$\mathscr M$, the vanishing of the first Pontryagin form $p_1$ has a 
topological significance. According to Hirzebruch's theorem
(\cite{Hirz}, \cite{Postnikov} pp 229-230), the Hirzebruch signature is
the integral of the first Pontryagin form; Gauss-Bonnet' theorem relates Euler's topological index
to a quadratic scalar of the curvature tensor \cite{Patterson}: 
$$   \tau (\mathscr M ) = \frac{1}{3}\int_{\mathscr M} p_1,  \qquad
\chi(\mathscr M)= \frac{3}{\pi^2}\int_{\mathscr M} (R^2-4R_{jik}R^{jk} +R_{jklm}R^{jklm}) $$
The two topological invariants coincide: $\tau = \chi$ mod.2 (\cite{Nakahara} page 465). 
\end{remark}
For (CQR)$_n$ manifolds the additional condition \eqref{eq2.13} has been 
investigated, and makes them conformally recurrent \cite{BuchRot}\cite{Prv90}:
\begin{equation}
\nabla_i C_{jklm} = 4 A_i C_{jklm}. \label{eq2.14} 
\end{equation}

\section{\bf (CQR)$_4$ Lorentzian manifolds}
We show that (CQR)$_4$ Lorentzian manifolds have the property
\eqref{eq2.13}. Then, besides being conformally harmonic
\eqref{eq2.2}, they are also conformally recurrent \eqref{eq2.14}, and
various geometric consequences follow.\\
Let us quote the useful general lemma:
\begin{lem} \label{lem3.1} {\rm (see \cite {LovRun} page 128)}\\
The Weyl tensor of a $n=4$ pseudo-Riemannian manifold satisfies 
the identity:
\begin{align}
 \delta_r^i C^{jk}{}_{ st} + \delta_t^i C^{jk}{}_{rs} + 
\delta_s^i C^{jk}{}_{tr} + \delta_r^k C^{ij}{}_{st} + 
\delta_t^k C^{ij}{}_{rs}  \nonumber \\
+\delta_s^k C^{ij}{}_{tr}
+\delta_r^j C^{ki}{}_{st} + \delta_t^j C^{ki}{}_{rs} + 
\delta_s^j C^{ki}{}_{tr} = 0 \label{eq3.1}
\end{align}
\end{lem}
In particular, in a (CQR)$_4$ manifold, on multiplying the equation by $A_j$ and using 
\eqref{eq2.1} we obtain eq.\eqref{eq2.13}:
$$ A_r C^{ki}{}_{st} + A_t C^{ki}{}_{rs} + A_s C^{ki}{}_{tr} = 0 $$ 
\begin{prop}
A (CQR)$_4$ pseudo-Riemannian manifold is conformally recurrent, 
the fundamental vector is null, and the first Pontryagin form vanishes.
\begin{proof}
Eq.\eqref{eq2.13} is verified, then \eqref{eq2.14} follows by property 
\eqref{eq1.1}. Contraction with $A^r$ gives that $A$ is null.
The first Pontryagin form is zero by prop. \ref{prop2.9}.
\end{proof}
\end{prop}

\noindent
Next we need another general lemma:
\begin{lem} \label{3.2} {\rm (\cite{DeFCl} page 46, \cite{Hall} page 148)}\\
In a 4-dimensional pseudo-Riemannian manifold let $A$ be a null vector and 
$B$ a vector orthogonal to $A$, $A_i B^i = 0$. Then $B$ is either space-like, 
or null and proportional to $A$, i.e. $B_j = \lambda A_j$ for some 
$\lambda \in R$.
\end{lem}
%
\noindent
We then prove:
\begin{prop} \label{th3.3}
In a 4-dimensional non conformally flat pseudo-Riemannian manifold, 
let $A$ and 
$B$ be vector fields such that $A_m C_{jkl}{}^m = 0$ and $B_m C_{jkl}{}^m = 0$.
Then $A^j A_j = 0$, $B_j B^j = 0$, and $B_j = \lambda A_j$ 
for some real $\lambda $ .
\begin{proof}
On multiplying eq.\eqref{eq3.1} by $A_j B^s$ we obtain
$(A_j B^j )C^{ki}{}_{tr} = 0$. Similarly, we obtain 
$( A^j A_j)C^{ki}{}_{tr} = 0$ and 
$( B^j B_j )C^{ki}{}_{tr} = 0$  (see \cite{LovRun} page 128). Then $A$ 
and $B$ are orthogonal null vectors. By lemma \ref{3.2} they are proportional.
\end{proof}
\end{prop}


The Bel-Debever version of Petrov's classification of Weyl tensors on 
4-dimensional Lorentzian manifolds (see \cite{Hall} page 196, \cite{Bel}, 
and \cite{Stephani}) is based on null vectors. 
By identifying it with the fundamental vector $A$, we 
may assert: 
\begin{prop}
On a non conformally flat (CQR)$_4$ Lorentzian manifold, the Weyl tensor is
type-N with respect to its fundamental vector. 
 \end{prop}
This shows that the assumption of a conformal motion in \cite{De2} is inessential.

\begin{prop}\label{th3.8}
On a non conformally flat (CQR)$_4$ Lorentzian manifold,\\
1) the fundamental vector $A$ is an eigenvector of the Ricci tensor;\\
2) the integral curves of the fundamental vector $A$ are geodesics;\\
3) if $A$ is closed, then $\nabla_i A^i =0$.
\begin{proof}
On multiplying \eqref{eq2.8} by $A^i$ and using eqs. \eqref{eq2.1}
one obtains $A^iR_{im} C_{jkl}{}^m =0$;  the equation
defines a vector $B_m = A^i R_{im}$ such that $B^m C_{jklm} =0$. 
By lemma 3.1 it is $A^i R_{im} = \lambda A_m$, i.e.  $A$ is an 
eigenvector of the Ricci tensor.\\
Next, multilpy eq.\eqref{eq2.3} by $A^i$ and obtain
$(A^i \nabla_i A^m )C_{jklm} = 0$. This equation and
 $A^mC_{jklm} = 0$ imply, by
the same arguments, $A^i (\nabla_i A_m) = \lambda A_m$. 
Therefore the integral curves of the vector A are geodesics 
(see \cite{Stephani} page 41).\\
Multiply \eqref{eq3.1} by $\nabla^p A^j$ and use \eqref{eq2.3}:
$$ (\nabla^p A_r) C^{ki}{}_{st} +(\nabla^p A_t) C^{ki}{}_{rs}
+(\nabla^p A_s) C^{ki}{}_{tr}=0. $$ 
Contraction of $s$ with $p$ and the closedness condition of $A$ give
 $(\nabla_s A^s ) C^{ki}{}_{tr}=0$.
\end{proof}
\end{prop}

From $A^mC_{jklm} = 0$ and $A^i R_{jm} = \lambda A_m$, a direct calculation
gives $A^mA^j R_{jklm} = (\lambda -\frac{1}{6}R)A_k A_l$
(see also Hall's theorem in \cite{Stephanietal} and \cite{WCT}). It follows 
that $$A_{[p} R_{k]jlm} A^mA^j =0$$  i.e. the Riemann's tensor is 
algebraically special (i.e. type II or D).


In their study of (PCS)$_n$ manifolds with proper conformal motion, 
De and Mazumdar \cite{De2} obtained the relations:
$(\nabla_i\sigma)(\nabla^i\sigma)=0$ and $C_{jkl}{}^m \nabla_m \sigma =0$.
In $n=4$ this proves that 
the space is type-N with respect to it and, by proposition \ref{th3.3} we have 
$\nabla_i \sigma = \lambda A_i$.  We conclude:

\begin{prop}
If a  non conformally flat (CQR)$_4$ Lorentzian manifold admits a proper 
conformal motion, then $\nabla_i\sigma=\lambda A_i$, i.e. the fundamental 
vector $A$ is closed.
\end{prop}

\end{document}